\documentclass[12pt]{amsart}
\usepackage{amsfonts}
\usepackage{amsmath,amssymb,amsthm}

\newcommand{\R}{\mathbb R}

\newcommand{\E}{\mathbb E}

\renewcommand{\span}{\mathrm{span}}
\newcommand{\tr}{\mathrm{tr}}

\newtheorem{thm}{Theorem}[section]

\theoremstyle{definition}

\theoremstyle{remark}

\newcommand{\ds}{\displaystyle}

\begin{document}
\title[MARGINALLY TRAPPED SURFACES WITH POINTWISE 1-TYPE
GAUSS MAP] {\bf MARGINALLY TRAPPED SURFACES
WITH POINTWISE 1-TYPE GAUSS MAP IN MINKOWSKI 4-SPACE}

\author{Velichka Milousheva}
\address{Bulgarian Academy of Sciences, Institute of Mathematics and Informatics,
Acad. G. Bonchev Str. bl. 8, 1113, Sofia, Bulgaria; "L. Karavelov"
Civil Engineering Higher School, 175 Suhodolska Str., 1373 Sofia,
Bulgaria} \email{vmil@math.bas.bg}



\maketitle

\begin{abstract}
A marginally trapped surface in the four-dimensional Minkowski
space is a spacelike surface whose mean curvature vector is
lightlike at each point. In the present paper we find all marginally trapped surfaces
with pointwise 1-type Gauss map. We prove that a marginally trapped surface is of
pointwise 1-type Gauss map if and only if it has parallel mean
curvature vector field.
\end{abstract}

\section{Introduction}

The concept of trapped surfaces was introduced by Roger Penrose
\cite{Pen}  and it plays an important role in General Relativity
for studying singularities and also for understanding the
evolution of black holes, the cosmic censorship hypothesis, the
Penrose inequality, etc. In Physics, a surface in a 4-dimensional
spacetime is called marginally trapped if it is closed, embedded,
spacelike and its mean curvature vector is lightlike at each point
of the surface. Recently, marginally trapped surfaces have been
studied from a mathematical viewpoint. In the mathematical
literature, it is customary to call a surface \emph{marginally
trapped} or \emph{quasi-minimal} \cite{Vran-Rosca, Chen-Garay} if
its mean curvature vector $H$ is lightlike at each point, and
removing the other hypotheses, i.e. the surface does not need to
be closed or embedded.

Classification results in 4-dimensional Lorentz space forms were
obtained imposing some extra conditions on the mean curvature
vector, the Gauss curvature or the second fundamental form. For
example, marginally trapped surfaces with positive relative
nullity in Lorenz space forms were classified by  B.-Y. Chen and
J. Van der Veken \cite{Chen-Veken-1}. They also proved the
non-existence of marginally trapped surfaces in Robertson-Walker
spaces with positive relative nullity \cite{Chen-Veken-2} and
classified marginally trapped surfaces with parallel mean
curvature vector in Lorenz space forms \cite{Chen-Veken-3}.

In \cite{Haesen-Ort-2}  S. Haesen and  M. Ortega classified
marginally trapped surfaces in Minkowski 4-space which are
invariant under spacelike rotations. Marginally trapped surfaces
in Minkowski 4-space which are invariant under boost
transformations  (hyperbolic rotations) were classified in
\cite{Haesen-Ort-1} and marginally trapped surfaces which are
invariant under the group of screw rotations (a group of Lorenz
rotations with an invariant lightlike direction) were studied  in
\cite{Haesen-Ort-3}.

In \cite{GM6} G. Ganchev and the present author developed the
invariant theory of marginally trapped surfaces in the
four-dimensional Minkowski space $\R^4_1$. Our approach to the
study of these surfaces was based on the principal lines generated
by the second fundamental form. Using the principal lines, we
introduced a geometrically determined moving frame field  at each
point of such a surface. The derivative  formulas for this frame
field imply the existence of
 seven invariant functions.  We proved that each marginally trapped surface is
determined up to a motion in $\R^4_1$ by these seven invariant
functions satisfying some natural conditions.

In the present paper we express  the Laplacian of the Gauss map of
a marginally trapped surface in terms of these invariant
functions. Imposing the condition that the surface has pointwise
1-type Gauss map, we obtain that three of the invariants are zero.
We give necessary and sufficient conditions for a marginally
trapped surface  to have pointwise 1-type Gauss map and find all
marginally trapped surfaces  with pointwise 1-type Gauss map. Our
main result states that a marginally trapped surface is of
pointwise 1-type Gauss map if and only if it has parallel mean
curvature vector field.

\section{Invariants of a marginally trapped surface} \label{S:Invariants}

Let  $\R^4_1$ be the Minkowski space endowed with the metric
$\langle , \rangle$ of signature $(3,1)$ and $Oe_1e_2e_3e_4$ be a
fixed orthonormal coordinate system in $\R^4_1$ such that $e_1^2 =
e_2^2 = e_3^2 = 1$, $e_4^2 = -1$.  The standard flat metric is
given in local coordinates by $dx_1^2 + dx_2^2 + dx_3^2 -dx_4^2.$

A surface $M^2$ in $\R^4_1$ is said to be \emph{spacelike} if
$\langle , \rangle$ induces  a Riemannian metric $g$ on $M^2$,
i.e. at each point $p$ of a spacelike surface $M^2$ we have the
following decomposition
$$\R^4_1 = T_pM^2 \oplus N_pM^2$$
with the property that the restriction of the metric
$\langle , \rangle$ onto the tangent space $T_pM^2$ is of
signature $(2,0)$, and the restriction of the metric $\langle ,
\rangle$ onto the normal space $N_pM^2$ is of signature $(1,1)$.

We denote by $\nabla'$ and $\nabla$ the Levi Civita connections on
$\R^4_1$ and $M^2$, respectively. Let $x$ and $y$ be vector fields
tangent to $M$ and let $\xi$ be a normal vector field. Then the
formulas of Gauss and Weingarten give  decompositions of the
vector fields $\nabla'_xy$ and $\nabla'_x \xi$ into  tangent and
normal components:
$$\begin{array}{l}
\vspace{2mm}
\nabla'_xy = \nabla_xy + \sigma(x,y);\\
\vspace{2mm} \nabla'_x \xi = - A_{\xi} x + D_x \xi,
\end{array}$$
which define the second fundamental tensor $\sigma$, the normal
connection $D$ and the shape operator $A_{\xi}$ with respect to
$\xi$. The  mean curvature vector  field $H$ of  $M^2$ is defined
as $H = \ds{\frac{1}{2}\,  \tr\, \sigma}$. Thus, if $M^2$ is a
spacelike surface and $\{x,y\}$ is a local orthonormal frame  of
the tangent bundle, the mean curvature vector  field is given by
the formula $H = \ds{\frac{1}{2} \left(\sigma(x,x) +
\sigma(y,y)\right)}$.

Let $M^2: z=z(u,v), \,\, (u,v) \in \mathcal{D}$ $(\mathcal{D}
\subset \R^2)$  be a local parametrization on a spacelike surface
in $\R^4_1$. The tangent space at an arbitrary point $p$ of $M^2$
is $T_pM^2 = \span \{z_u,z_v\}$. Since $M^2$ is spacelike,
$\langle z_u,z_u \rangle > 0$, $\langle z_v,z_v \rangle > 0$. We
use the standard denotations $E(u,v)=\langle z_u,z_u \rangle, \;
F(u,v)=\langle z_u,z_v \rangle, \; G(u,v)=\langle z_v,z_v \rangle$
for the coefficients of the first fundamental form
$$I(\lambda,\mu)= E \lambda^2 + 2F \lambda \mu + G \mu^2,\quad
\lambda, \mu \in \R$$ and set $W=\sqrt{EG-F^2}$. Let us choose a
normal frame field $\{n_1, n_2\}$ such that $\langle n_1, n_1
\rangle =1$, $\langle n_2, n_2 \rangle = -1$, and the quadruple
$\{z_u,z_v, n_1, n_2\}$ is positively oriented in $\R^4_1$. Then
we have the following derivative formulas:
$$\begin{array}{l}
\vspace{2mm} \nabla'_{z_u}z_u=z_{uu} = \Gamma_{11}^1 \, z_u +
\Gamma_{11}^2 \, z_v + c_{11}^1\, n_1 - c_{11}^2\, n_2;\\
\vspace{2mm} \nabla'_{z_u}z_v=z_{uv} = \Gamma_{12}^1 \, z_u +
\Gamma_{12}^2 \, z_v + c_{12}^1\, n_1 - c_{12}^2\, n_2;\\
\vspace{2mm} \nabla'_{z_v}z_v=z_{vv} = \Gamma_{22}^1 \, z_u +
\Gamma_{22}^2 \, z_v + c_{22}^1\, n_1 - c_{22}^2\, n_2,\\
\end{array}$$
where $\Gamma_{ij}^k$ are the Christoffel's symbols and the functions $c_{ij}^k, \,\, i,j,k = 1,2$ are given by
$$\begin{array}{lll}
\vspace{2mm}
c_{11}^1 = \langle z_{uu}, n_1 \rangle; & \qquad   c_{12}^1 = \langle z_{uv}, n_1 \rangle; & \qquad  c_{22}^1 = \langle z_{vv}, n_1 \rangle;\\
\vspace{2mm}
c_{11}^2 = \langle z_{uu}, n_2 \rangle; & \qquad  c_{12}^2 = \langle z_{uv}, n_2 \rangle; & \qquad
c_{22}^2 = \langle z_{vv}, n_2 \rangle.
\end{array} $$

Obviously, the surface $M^2$ lies in a 2-plane if and only if
$M^2$ is totally geodesic, i.e. $c_{ij}^k=0, \; i,j,k = 1, 2.$ So,
we assume that at least one of the coefficients $c_{ij}^k$ is not
zero.

Let $X=\lambda z_u+\mu z_v, \,\, (\lambda,\mu)\neq(0,0)$ be a
tangent vector at a point $p \in M^2$. The second fundamental form
$II$ of the surface $M^2$ at the point $p$ is introduced by the
formula
$$II(\lambda,\mu)=L\lambda^2+2M\lambda\mu+N\mu^2,$$
where the functions $L$, $M$, and $N$ are defined as follows:
\begin{equation}\notag \label{Eq-1}
L = \ds{\frac{2}{W}} \left|%
\begin{array}{cc}
\vspace{2mm}
  c_{11}^1 & c_{12}^1 \\
  c_{11}^2 & c_{12}^2 \\
\end{array}%
\right|; \quad
M = \ds{\frac{1}{W}} \left|%
\begin{array}{cc}
\vspace{2mm}
  c_{11}^1 & c_{22}^1 \\
  c_{11}^2 & c_{22}^2 \\
\end{array}%
\right|; \quad
N = \ds{\frac{2}{W}} \left|%
\begin{array}{cc}
\vspace{2mm}
  c_{12}^1 & c_{22}^1 \\
  c_{12}^2 & c_{22}^2 \\
\end{array}%
\right|.
\end{equation}
The second fundamental form $II$ is invariant up to the
orientation of the tangent space or the normal space of the
surface.

The condition $L = M = N = 0$  characterizes points at which the
space $\{\sigma(x,y):  x, y \in T_pM^2\}$ is one-dimensional. We
call such points  \emph{flat points} of the surface \cite{GM4,
GM5}. These points are analogous to flat points in the theory of
surfaces in $\R^3$. In \cite{Lane} and \cite{Little} such points
are called inflection points. The notion of an inflection point is
introduced for 2-dimensional surfaces in a 4-dimensional affine
space $\mathbb{A}^4$. E. Lane \cite{Lane} has shown that every
point of a surface in $\mathbb{A}^4$ is an inflection point if and
only if the surface is developable or lies in a 3-dimensional
space. Further we consider surfaces free of flat points, i.e. $(L,
M, N) \neq (0,0,0)$.

The second fundamental form $II$ determines conjugate, asymptotic,
and principal tangents at a point $p$ of $M^2$ in the standard
way. A line $c: u=u(q), \; v=v(q); \; q\in J \subset \R$ on $M^2$
is said to be  a \textit{principal line}, if its tangent at any
point is principal.

It is interesting to note that the ''\emph{umbilical}'' points,
i.e. points at which the coefficients of the first and the second
fundamental forms are proportional ($L = \rho E, \, M = \rho F,
\,N = \rho G, \, \rho \neq 0$), are exactly the points at which
the mean curvature vector $H$ is zero. So, the spacelike surfaces
consisting of ''umbilical'' points in $\R^4_1$ are exactly the
 surfaces with zero mean curvature. If $M^2$ is a spacelike surface free of ''umbilical''
 points ($H \neq 0$ at each point), then there exist exactly two principal tangents.

\vskip 2mm Now, let $M^2: z=z(u,v), \,\, (u,v) \in \mathcal{D}$ be
a marginally trapped surface. Then the mean curvature vector is
lightlike at each point of the surface, i.e. $\langle H, H \rangle
= 0$, $H \neq 0$. Hence there exists a pseudo-orthonormal normal
frame field $\{n_1, n_2\}$, such that $n_1 = H$ and
$$\langle n_1, n_1 \rangle = 0; \quad \langle n_2, n_2 \rangle = 0; \quad \langle n_1, n_2 \rangle = -1.$$
We assume that $M^2$  is  free of flat points, i.e. $(L,M,N) \neq
(0,0,0)$. Then at each point of the surface there exist principal
lines and without loss of generality we assume that $M^2$ is
parameterized by principal lines. Let us denote $x =
\ds{\frac{z_u}{\sqrt{E}}}$, $y = \ds{\frac{z_v}{\sqrt{G}}}$. Thus
we obtain a special frame field $\{x,y,n_1,n_2\}$ at each point $p
\in M^2$, such that $x, y$ are unit spacelike vector fields
collinear with the principal directions; $n_1, n_2$ are lightlike
vector fields, $\langle n_1,\, n_2 \rangle = -1$, and $n_1$ is the
mean curvature vector field. We call this frame field a
\emph{geometric frame field} of $M^2$.

With respect to the geometric frame field we have the following
derivative formulas of $M^2$:
\begin{equation} \label{Eq-5}
\begin{array}{ll}
\vspace{2mm} \nabla'_xx=\quad \quad \quad
\gamma_1\,y+\,(1+\nu)\,n_1; & \qquad
\nabla'_x n_1= \quad\quad \quad \quad  \mu\,y+\beta_1\,n_1;\\
\vspace{2mm} \nabla'_xy=-\gamma_1\,x\quad \quad \; + \,\quad
\lambda\,n_1 \; + \mu\,n_2;  & \qquad
\nabla'_y n_1=\mu \,x \quad\quad \quad \quad +\beta_2\,n_1;\\
\vspace{2mm} \nabla'_yx=\quad\quad \;\, -\gamma_2\,y  + \quad
\lambda\,n_1 \; +\mu\,n_2;  & \qquad
\nabla'_xn_2= (1+\nu) \, x + \lambda \,y \quad \quad  -\beta_1\,n_2;\\
\vspace{2mm} \nabla'_yy=\;\;\gamma_2\,x \quad\quad\;\;\,
+(1-\nu)\,n_1; & \qquad \nabla'_y n_2= \lambda \, x  + (1-\nu) \,y
 \quad \quad -\beta_2\,n_2,
\end{array}
\end{equation}
where $\nu = - \langle \ds{\frac{\sigma(x,x) - \sigma(y,y)}{2}, n_2} \rangle$,
$\lambda =-  \langle \sigma(x,y), n_2\rangle$, $\mu = -\langle \sigma(x,y), n_1\rangle$,
 $\gamma_1 = - y(\ln \sqrt{E})$, $\gamma_2 = - x(\ln \sqrt{G})$, $\beta_1 = - \langle \nabla'_x n_1, n_2\rangle$,
$\beta_2 = - \langle \nabla'_y n_1, n_2\rangle$.

Note that the functions $\nu,\,\lambda, \, \mu,\, \gamma_1,\, \gamma_2,\,
\beta_1, \, \beta_2$ are invariants of
the surface determined by the principal directions.
In \cite{GM6} we proved the fundamental theorem for marginally trapped
 surfaces in $\R^4_1$, which states that each marginally trapped surface free of flat points is
determined up to a motion in $\R^4_1$ by these seven invariant
functions satisfying some natural conditions.

Formulas \eqref{Eq-5} imply that the Gauss curvature $K$ and normal curvature  $\varkappa$ of $M^2$
are expressed by the functions $\nu, \lambda$, and $\mu$ as
follows:

\begin{equation} \label{Eq-3}
 \qquad K = 2 \lambda \mu; \qquad \varkappa = - 2 \mu \nu.
\end{equation}

\section{Surfaces  with pointwise 1-type Gauss map}\label{S:Pointwise}

An isometric immersion $x:M$ $\rightarrow $ $\mathbb{E}^{m}$ of a
submanifold $M$ in the Euclidean space $\mathbb{E}^{m}$ or
pseudo-Euclidean space $\E^m_s$ with signature $(s,m-s)$ is said
to be of \emph{finite type}  \cite{Ch1}, if $x$ identified with
the position vector field of $M$ in $\mathbb{E}^{m}$ or $\E^m_s$
can be expressed as a finite sum of eigenvectors of the Laplacian
$\Delta $ of $M$, i.e.
\begin{equation*}
x=x_{0}+\sum_{i=1}^{k}x_{i},
\end{equation*}
where $x_{0}$ is a constant map, $x_{1},x_{2},...,x_{k}$ are non-constant maps
such that $\Delta x_i=\lambda _{i}x_{i},$ $\lambda _{i}\in \mathbb{R}$, $1\leq
i\leq k$. If $\lambda _{1},\lambda _{2},...,\lambda _{k}$ are different,
then $M$ is said to be of \emph{$k$-type}.

The notion of finite type immersion is naturally extended to the
Gauss map $G$ on $M$ by  B.-Y. Chen and  P. Piccinni \cite{CP}.
Thus,  a submanifold $M$ of an  Euclidean or pseudo-Euclidean
space has \emph{1-type Gauss map} $G$, if $G$ satisfies $\Delta
G=a (G+C)$ for some $a \in \mathbb{R}$ and some constant vector
$C$.

However, the Laplacian of the Gauss map of some typical well-known
surfaces in the three-dimensional Euclidean space $\mathbb{E}^{3}$
such as the helicoid, the catenoid and the right cone takes a
somewhat different form, namely, $\Delta G=\phi (G+C)$ for some
non-constant function $\phi$ and some constant vector $C$. It
looks like an eigenvalue problem, but the function $\phi$ turns
out to be non-constant. Therefore, it is worth studying the class
of surfaces satisfying such an equation.

We use the following definition: a submanifold $M$ of the
 Euclidean space $\mathbb{E}^{m}$ or
pseudo-Euclidean space $\E^m_s$ is said to have \emph{pointwise
1-type Gauss map} if its Gauss map $G$ satisfies
\begin{equation}\notag
\Delta G=\phi (G+C)  \label{A1}
\end{equation}
for some  smooth function $\phi$ on $M$ and a constant vector $C$.
 A pointwise 1-type Gauss map is called \emph{proper} if the function $\phi $
 is non-constant. A submanifold with pointwise 1-type
Gauss map is said to be of \emph{the first kind} if the vector $C$ is
zero. Otherwise, the pointwise 1-type Gauss map is said to be o\emph{f the
second kind}.

In \cite{CK} M. Choi and Y.  Kim characterized the helicoid in
terms of pointwise 1-type Gauss map of the first kind. Together
with B.-Y. Chen they  proved that the class of surfaces of
revolution with pointwise 1-type Gauss map of the first kind
coincides with the class of surfaces of revolution with constant
mean curvature and characterized the rational surfaces of
revolution with pointwise 1-type Gauss map  \cite{CCK}.

Tensor product surfaces with pointwise 1-type Gauss map and
Vranceanu rotational surfaces with pointwise 1-type Gauss map in
the four-dimensional Euclidean space $\E^4$ were studied in
\cite{A2} and \cite{A1}, respectively.

In \cite{KY3} Y.  Kim and  D. Yoon studied ruled surfaces with
1-type Gauss map in Minkowski space $\E^m_1$ and gave a complete
classification of null scrolls with 1-type Gauss map. The
classification of ruled surfaces with pointwise 1-type Gauss map
of the first kind in Minkowski space $\E^3_1$ was given in
\cite{KY2}. Ruled surfaces with pointwise 1-type Gauss map  of the
second kind in Minkowski 3-space were classified in \cite{CKY}.

\vskip 1mm Recall that the Gauss map $G$ of a submanifold $M$ of
$\E^m$ is defined as follows. Let $G(n,m)$ be the Grassmannian
manifold consisting of all oriented $n$-planes through the origin
of $\mathbb{E}^{m}$ and $\wedge ^{n}\mathbb{E}^{m}$ be the vector
space obtained by the exterior product of $n$ vectors in
$\mathbb{E}^{m}$. In a natural way, we can identify $\wedge
^{n}\mathbb{E}^{m}$ with the Euclidean space $\mathbb{E}^{N}$,
where $N=\left(
\begin{array}{c}
m \\
n
\end{array}
\right)$.
 Let $\left\{ e_{1},...,e_{n},e_{n+1},\dots,e_{m}\right\} $ be a
 local orthonormal frame field in $\mathbb{E}^{m}$ such that $e_{1},e_{2},\dots,$ $e_{n}$ are tangent to $M$ and
 $e_{n+1},e_{n+2},\dots,e_{m}$ are
 normal to $M$.
The map $G:M\rightarrow G(n,m)$ defined by $%
G(p)=(e_{1}\wedge e_{2}\wedge \dots \wedge $ $e_{n})(p)$ is called the \emph{Gauss
map} of $M$. It is a smooth map which carries a point $p$ in $M$ into the
oriented $n$-plane in $\mathbb{E}^{m}$ obtained by the parallel translation
of the tangent space of $M$ at $p$ in $\mathbb{E}^{m}$.

In a similar way one can consider the Gauss map of  a submanifold $M$ of pseudo-Euclidean space  $\E^m_s$.

For any  function $f$ on $M$ the Laplacian of $f$ is given by the formula

\begin{equation}\notag
\Delta f =-\sum_{i}(\nabla'_{e_{i}}\nabla'
_{e_{i}}f -\nabla'_{\nabla _{e_{i}}e_{i}}f ),
\end{equation}
where $\nabla'$ is the Levi-Civita connection of $\E^m$ or $\E^m_s$  and $\nabla$ is the induced connection on $M$.

\section{Main result}

In this section we shall study marginally trapped surfaces with pointwise 1-type Gauss map.

Let $M^2$ be a marginally trapped surface free of flat points and
$\{x,y,n_1,n_2\}$ be the geometric frame field of $M^2$, defined
in Section \ref{S:Invariants}. The geometric frame field
$\{x,y,n_1,n_2\}$ generates the following frame of the
Grassmannian manifold: $$\{x \wedge y, x \wedge n_1, x \wedge n_2,
y \wedge n_1, y \wedge n_2, n_1 \wedge n_2\}.$$ The indefinite
inner product on the Grassmannian manifold is given by
\begin{equation}\notag
\langle e_{i_1} \wedge e_{i_2}, f_{j_1} \wedge f_{j_2}  \rangle =
\det \left( \langle e_{i_k}, f_{j_l}  \rangle  \right).
\end{equation}

The Gauss map $G$ of $M^2$ is defined by $G(p)=(x\wedge y)(p)$, $p
\in M^2$. Then the Laplacian of the Gauss map is given by the
formula
\begin{equation}\label{Eq-7}
\Delta G = - \nabla'_x\nabla'_x G + \nabla'_{\nabla_x x} G -
\nabla'_y\nabla'_y G + \nabla'_{\nabla_y y} G.
\end{equation}

The derivative formulas \eqref{Eq-5} of $M^2$ imply the following
equalities for the invariants $\nu,\,\lambda, \, \mu,\,
\gamma_1,\, \gamma_2,\, \beta_1, \, \beta_2$ of the surface:
\begin{equation} \label{Eq-6}
\begin{array}{l}
\vspace{2mm}
 x(\mu) - 2\mu\, \gamma_2  - \mu\,\beta_1 =0;\\
\vspace{2mm}
y(\mu) - 2\mu\, \gamma_1 - \mu\,\beta_2 = 0;\\
\vspace{2mm}
x(\lambda) - y(\nu) - 2\lambda\, \gamma_2 + 2\nu\, \gamma_1 + \lambda\,\beta_1 - (1+ \nu)\,\beta_2 = 0;\\
\vspace{2mm}
x(\nu) + y(\lambda) - 2\lambda\, \gamma_1 - 2\nu\, \gamma_2 - (1- \nu)\,\beta_1 + \lambda\,\beta_2 = 0;\\
\vspace{2mm}
x(\beta_2) - y(\beta_1) + 2 \nu\,\mu + \gamma_1\,\beta_1 - \gamma_2\,\beta_2 = 0.
\end{array}
\end{equation}

Using equalities \eqref{Eq-5} and \eqref{Eq-7} we calculate the Laplacian of
the Gauss map:
\begin{equation}\label{Eq-8}
\begin{array}{ll}
\vspace{2mm}\Delta G =& -4 \lambda \,\mu\, x \wedge y - \left(x(\lambda) - y(\nu) - 2\lambda\, \gamma_2 + 2\nu\, \gamma_1 + \lambda\,\beta_1 + (1- \nu)\,\beta_2 \right)\, x \wedge n_1 \\
\vspace{2mm}
 & - \left(  x(\mu) - 2\mu\, \gamma_2  - \mu\,\beta_1 \right)\, x \wedge n_2 +  \left( y(\mu) - 2\mu\, \gamma_1 - \mu\,\beta_2  \right)\, y \wedge n_2\\
 \vspace{2mm}
 & + \left(x(\nu) + y(\lambda) - 2\lambda\, \gamma_1 - 2\nu\, \gamma_2 + (1+ \nu)\,\beta_1 + \lambda\,\beta_2 \right)\, y \wedge n_1 \\
\vspace{2mm}
 & - 4 \mu\, \nu\, n_1 \wedge n_2.
 \end{array}
 \end{equation}

Equalities \eqref{Eq-6} and \eqref{Eq-8} imply that the  Laplacian of
the Gauss map of a marginally trapped surface free of flat points is expressed in terms of the invariants
of the surface by the following formula:
\begin{equation}\label{Eq-9}
\Delta G = -4 \lambda \,\mu\, x \wedge y - 2 \beta_2\, x \wedge n_1 + 2 \beta_1\, y \wedge n_1 - 4 \mu\, \nu\, n_1 \wedge n_2.
\end{equation}

Using \eqref{Eq-3} we can rewrite \eqref{Eq-9} in the form:
\begin{equation}\label{Eq-10}\notag
\Delta G = - 2K \, x \wedge y - 2 \beta_2\, x \wedge n_1 + 2 \beta_1\, y \wedge n_1 + 2 \varkappa\, n_1 \wedge n_2,
\end{equation}
where $K$ and $\varkappa$ are the Gauss curvature and the normal curvature, respectively.

\vskip 2mm
Now we shall find all marginally trapped surfaces with pointwise 1-type Gauss map.

\begin{thm}
Let $M^2$ be a marginally trapped surface free of flat points. Then
$M^2$ is of pointwise 1-type Gauss map if and
only if $M^2$ has parallel mean curvature vector field.
\end{thm}

\vskip 2mm \noindent \emph{Proof:}
Let $M^2$ be a marginally trapped surface free of flat points. Then the  Laplacian of
the Gauss map is  expressed  by  formula \eqref{Eq-9}.

In the case when the Gauss curvature is non-zero at a point $p\in
M^2$, we have that $\lambda \neq 0$ in a neighbourhood of $p$.
Then  the Laplacian of the Gauss map can be written as
\begin{equation}\label{Eq-11}
\Delta G = - 4\lambda \mu \, G - 4\lambda \mu \, T,
\end{equation}
where $T = \ds{\frac{\beta_2}{2\lambda \mu}\, x \wedge n_1 - \frac{\beta_1}{2\lambda \mu}\, y \wedge n_1 +\frac{ \nu}{\lambda}\, n_1 \wedge n_2}$.

The condition that the
surface has  pointwise 1-type Gauss map is
$$\Delta G=\phi (G+C)$$
for some  smooth function $\phi$ on $M^2$ and a constant vector $C$.
Using \eqref{Eq-11} we obtain that $M^2$ is of  pointwise 1-type Gauss map if and only if $T = const$.

By the use of formulas \eqref{Eq-5} we calculate that
\begin{equation}\label{Eq-12}
\begin{array}{ll}
\vspace{2mm}
\nabla'_x T =& \ds{\frac{\beta_2}{2 \lambda}\, x \wedge y +\left(x(\frac{\beta_2}{2 \lambda \mu}) + \frac{\beta_1 \beta_2 + \beta_1 \gamma_1}{2 \lambda \mu} - \frac{\nu (1+\nu)}{ \lambda} \right)\, x \wedge n_1 }\\
\vspace{2mm}
 & + \ds{ \left( - x(\frac{\beta_1}{2 \lambda \mu}) + \frac{\beta_2 \gamma_1 - (\beta_1)^2}{2 \lambda \mu} - \nu \right)\,  y \wedge n_1 +  \frac{\nu\mu}{\lambda}\, y \wedge n_2}\\
 \vspace{2mm}
  & + \ds{\left(\frac{\beta_1}{2\lambda} + x(\frac{\nu}{\lambda}) \right)\, n_1 \wedge n_2};\\
 \vspace{2mm}
\nabla'_y T =& \ds{\frac{\beta_1}{2 \lambda}\, x \wedge y +\left(y(\frac{\beta_2}{2 \lambda \mu}) + \frac{(\beta_2)^2 - \beta_1 \gamma_2}{2 \lambda \mu} - \nu \right)\, x \wedge n_1 }\\
\vspace{2mm}
 & + \ds{\frac{\nu\mu}{\lambda}\, x \wedge n_2 + \left( - y(\frac{\beta_1}{2 \lambda \mu}) - \frac{\beta_1 \beta_2 + \beta_2 \gamma_2}{2 \lambda \mu} - \frac{\nu(1-\nu)}{\lambda} \right)\,  y \wedge n_1 }\\
 \vspace{2mm}
  & + \ds{\left(y(\frac{\nu}{\lambda})  - \frac{\beta_2}{2\lambda} \right)\, n_1 \wedge n_2}.
\end{array}
 \end{equation}

Equalities \eqref{Eq-12} imply that $T = const$ if and only if $\beta_1 =0$, $\beta_2 =0$, $\nu = 0$ in this neighbourhood.

\vskip 2mm

If $M^2$ is flat, i.e. $\lambda=0$ at each point, then the
Laplacian of the Gauss map is
\begin{equation}\notag
\Delta G = -2 T_0,
\end{equation}
where $T_0 = \ds{\beta_2\, x \wedge n_1 - \beta_1\, y \wedge n_1 + 2\mu \nu\, n_1 \wedge n_2}$.
Using \eqref{Eq-5} we find
\begin{equation}\notag
\begin{array}{ll}
\vspace{2mm}
\nabla'_x T_0 =& \ds{\mu \beta_2\, x \wedge y +\left(x(\beta_2) + \beta_1 \beta_2 + \beta_1 \gamma_1 - 2\mu\nu (1+\nu) \right)\, x \wedge n_1 }\\
\vspace{2mm}
 & + \ds{ \left( - x(\beta_1) + \beta_2 \gamma_1 - (\beta_1)^2\right)\,  y \wedge n_1 + 2\mu^2 \nu\, y \wedge n_2}\\
 \vspace{2mm}
  & + \ds{\left(x(2\mu\nu) + \mu \beta_1 \right)\, n_1 \wedge n_2};\\
\vspace{2mm}
\nabla'_y T_0 =& \ds{\mu \beta_1\, x \wedge y +\left(y(\beta_2) + (\beta_2)^2  - \beta_1\gamma_2 \right)\, x \wedge n_1 }\\
\vspace{2mm}
 & + \ds{2\mu^2 \nu \, x \wedge n_2 + \left( - y(\beta_1) -\beta_1 \beta_2 - \beta_2 \gamma_2 - 2\mu\nu(1-\nu) \right)\,  y \wedge n_1 }\\
 \vspace{2mm}
  & + \ds{\left(y(2\mu\nu)  -\mu \beta_2 \right)\, n_1 \wedge n_2},
\end{array}
 \end{equation}
which imply again that $T_0 = const$ if and only if $\beta_1 =0$, $\beta_2 =0$, $\nu = 0$.

\vskip 1mm Marginally trapped surfaces satisfying $\beta_1 =
\beta_2 =0$ have parallel mean curvature vector field, since $DH
=0$ holds identically in view of  \eqref{Eq-5}. From \eqref{Eq-3}
it follows that the  equality $\nu = 0$ is equivalent to
$\varkappa = 0$, i.e. the surface is of flat normal connection.
The last equality of \eqref{Eq-6} implies that all marginally
trapped surfaces with parallel mean curvature vector field
($\beta_1 = \beta_2 =0$) have flat normal connection ($\nu = 0$).

\vskip 1mm Finally we obtain that $M^2$ is of pointwise 1-type
Gauss map if and only if $M^2$ has parallel mean curvature vector
field. \qed

\vskip 3mm
Note that the Laplacian of the Gauss map of each marginally trapped surface $M^2$ with pointwise 1-type Gauss map
is expressed as follows:
$$\Delta G = -4 \lambda \,\mu\,G = -2 K \,G,$$
where $K$ is the Gauss curvature of $M^2$.
Hence,  $M^2$ is of the first kind ($C =0$). Moreover, $M^2$ is proper
($\phi \neq const$) in the case of non-constant Gauss curvature.

\vskip 3mm

The class of marginally trapped surfaces with parallel mean
curvature vector field, was classified by  B.-Y. Chen and J. Van
Der Veken in \cite{Chen-Veken-3}. Combining their classification
with our result, we obtain all marginally trapped surfaces with
pointwise 1-type Gauss map.

\end{document}